\documentclass{amsart}

\DeclareMathOperator{\val}{val}

\DeclareMathOperator{\dom}{dom}
 
\DeclareMathOperator{\id}{id}

\DeclareMathOperator{\ad}{ad}
\DeclareMathOperator{\defect}{def}

\DeclareMathOperator{\ssc}{sc}
\DeclareMathOperator{\fr}{fr} 
\DeclareMathOperator{\rk}{rk}
\DeclareMathOperator{\Int}{Int}

\usepackage{amssymb}

\numberwithin{equation}{subsection}

\newtheorem{theorem}{Theorem}[subsection]
\newtheorem{corollary}[theorem]{Corollary}
\newtheorem{lemma}[theorem]{Lemma}
\newtheorem{proposition}[theorem]{Proposition}

\theoremstyle{definition}

\begin{document}
\title[Dimensions of Newton strata]{Dimensions of Newton strata in the
adjoint quotient of reductive groups}

\author{Robert E. Kottwitz}
\address{Department of Mathematics\\ University of Chicago\\ 5734 University
Avenue\\ Chicago, Illinois 60637}
\email{kottwitz@math.uchicago.edu}
\thanks{Partially supported by NSF Grant DMS-0245639}

\dedicatory{Dedicated to R.~MacPherson on  his 60th birthday}

\maketitle

\section{Introduction} 
In \cite{chai00} Chai made a conjecture on the codimensions
of Newton strata in Shimura varieties, which then led Rapoport 
\cite{rapoport05} to his conjecture on dimensions
of affine Deligne-Lusztig varieties inside affine Grassmannians.  The main
goal of this paper is to show that exactly the same codimensions arise in a
simpler context, that of the Newton stratification in the adjoint quotient
of a reductive group. Along the way we study this stratification and then
introduce the notion of
\emph{defect}, which we use to rewrite the codimension formula without
having to use the greatest integer function. 

The influence of the work of Chai and Rapoport on this paper is obvious.
Less obvious is the influence of some joint work (as yet unpublished) with
M.~Goresky and R.~MacPherson on codimensions of root-valuation strata in
Lie algebras over Laurent power series fields. It is an especially great 
pleasure to acknowledge the influence of Goresky and MacPherson in a paper
like this one, dedicated to Bob MacPherson on his 60th birthday. In
addition I would like to thank T.~Haines and M.~Sabitova for some very 
helpful comments on an earlier version of this paper.

\subsection{Preliminaries concerning the discretely valued field $F$}

Let $\mathcal O$ be a complete discrete valuation ring  with uniformizing
element $\pi$ and algebraically closed residue field 
$k:=\mathcal O/\pi\mathcal O$. 
In this paper it turns out to be more convenient to use the negative of the
usual valuation. Thus our valuation on the fraction field $F$ of $\mathcal
O$ is a surjective homomorphism 
\[
\val:F^\times \twoheadrightarrow \mathbb Z,
\]
normalized so that $\val(\pi)=-1$. It is then appropriate to put
$\val(0)=-\infty$. With these conventions the ultrametric inequality becomes
\[
\val(a+b) \le \max\{\val(a),\val(b)\}
\]
with equality when the valuations of $a$ and $b$ are distinct.

 Let $\bar F$ be an
algebraic closure of $F$ and extend the given valuation on $F$ to a
valuation on $\bar F$. Thus we now have a surjective
homomorphism $\val:\bar F^\times \twoheadrightarrow \mathbb Q$.  

\subsection{Classical Newton polygons}
Let $f(T)=T^n + a_1T^{n-1} + \dots +a_{n-1}T +a_n$ be a monic polynomial of
degree
$n$ with coefficients in $F$ and non-zero constant term $a_n$. The classical
theory of Newton polygons (see
\cite{robba80} for example) concerns the relation between the valuations of
the
$n$ roots of
$f(T)$ in $\bar F$ and the valuations of the coefficients of $f(T)$. Recall
that the valuations of the roots of $f$ are determined by the valuations of
the coefficients of
$f$ (but not the other way around). We now review the theory in greater
detail. 

Since $a_n\ne0$, all $n$ roots of $f(T)$ are non-zero, and their valuations
$\nu_1,\dots,\nu_n$, when arranged in decreasing order $\nu_1 \ge \dots
\ge \nu_n$, give us a well-defined point $\nu=(\nu_1,\dots,\nu_n) \in
\mathbb Q^n$.

We get another $n$-tuple $d=(d_1,\dots,d_n)$ by putting $d_i=\val(a_i)$. Thus
$d$ lies in $\tilde{\mathbb Z}^{n-1} \times \mathbb Z$, where $\tilde{\mathbb
Z}:=\mathbb Z \cup \{-\infty\}$. Up to a sign, $a_i$ is the $i$-th
elementary symmetric function of the roots, and it is clear that the
valuation of each monomial contributing to $a_i$ has valuation $\le \nu_1 +
\nu_2+\dots +\nu_i$ with equality occurring for at least one monomial.
Therefore 
\begin{equation}\label{eq.c.n.ineq}
d_i \le\nu_1+\dots + \nu_i \text{ \quad for $i=1,\dots,n-1$}
\end{equation}
and 
\begin{equation}\label{eq.c.n.eq}
d_i =\nu_1+\dots + \nu_i \text{ \quad for $i=n$.}
\end{equation}
Moreover, when $\nu_i > \nu_{i+1}$, exactly one monomial has 
valuation $\nu_1+\dots +\nu_i$, the rest having strictly smaller
valuations, so that equality holds in \eqref{eq.c.n.ineq} in this case. 

Thus,
as mentioned before, $\nu$ does not determine $d$. However, $d$ does 
determine $\nu$. For this we need the following
partial order on $\mathbb Q^n$: for $\mu,\mu' \in \mathbb Q^n$ we write
$\mu' 
\le \mu$ if and only if 
\begin{equation}\label{eq.root.ineq}
\mu'_1+\dots + \mu'_i \le \mu_1+\dots + \mu_i\text{ \quad for
$i=1,\dots,n-1$}
\end{equation}
and 
\begin{equation}\label{eq.root.eq}
\mu'_1+\dots + \mu'_i = \mu_1+\dots + \mu_i \text{ \quad for $i=n$}.
\end{equation}
It then turns out that $\nu$ is the unique smallest element in the set of
all $\mu \in \mathbb Q^n$ satisfying the three conditions
\begin{equation}
\mu_1 \ge \dots \ge \mu_n
\end{equation}
\begin{equation}\label{eq.c.n.ineq**}
d_i \le\mu_1+\dots + \mu_i \text{ \quad for $i=1,\dots,n-1$}
\end{equation} 
\begin{equation}\label{eq.c.n.eq**}
d_i =\mu_1+\dots + \mu_i \text{ \quad for $i=n$}.
\end{equation}
In classical language $\nu_1,\dots,\nu_n$ are the slopes of the Newton
polygon obtained from  the $n$-tuple $d$.

This story is related to the reductive group $GL_n$ in two ways. First,  an
element $g \in GL_n(F)$ has a characteristic polynomial
$f(T)$, whose roots are the $n$ eigenvalues of $g$. This provides  a natural
source of the kind of polynomials we have been considering. Second, as is
well-known, the partial order $\le$ and inequalities such as $\nu_1 \ge
\dots \ge\nu_n$ and $\nu_i>\nu_{i+1}$ occur most naturally in the context of
the root system for $GL_n$. 
Our first goal in this paper is to extend the
theory of Newton polygons to other reductive groups, as has already been
done for isocrystals with $G$-structure
\cite{kottwitz85},\cite{rapoport-richartz96},\cite{chai00}. This requires 
some preparation involving root systems and reductive groups, all of which
is standard, but which we review for the sake of completeness.

\subsection{Root theoretic preliminaries}\label{subsec.rtp}
Consider a split connected reductive group $G$ over
$\mathcal O$. We fix a split maximal torus $A$ in $G$ as well as a Borel
subgroup
$B=AU$ containing $A$, with $U$ denoting the unipotent radical of $B$. We
assume that $B$, $A$, $U$ are defined over $\mathcal O$, but most of the
time we will really be thinking about $G$, $A$, $B$, $U$ as algebraic groups
over $F$.  We also assume that the derived group of
$G$ is simply connected. We denote this derived group by $G_{\ssc}$. Thus
there is an exact sequence
\[
1 \to G_{\ssc} \to G \to D \to 1,
\]
where $D:=G/G_{\ssc}$, a split torus. We write $A_{\ssc}$ for the 
intersection of $A$ with $G_{\ssc}$. Thus $A_{\ssc}$ is a split maximal torus
in  $G_{\ssc}$, and there is an exact sequence
\[
1 \to A_{\ssc} \to A \to D \to 1,
\] 
giving rise to the exact sequence 
\[
0 \to X^*(D) \to X^*(A) \to X^*(A_{\ssc}) \to 0
\]
of character groups. 

Let $\alpha_1,\dots,\alpha_l \in X^*(A)$ denote the
simple roots of $A$, and let 
$
\alpha^\vee_1,\dots,\alpha^\vee_l \in
X_*(A_{\ssc}) 
$
denote the corresponding simple coroots. We write $\Delta$ for the set of
simple roots. The simple coroots form a
$\mathbb Z$-basis for $X_*(A_{\ssc})$, and we write
$\varpi_1,\dots,\varpi_l$ for the corresponding dual basis for
$X^*(A_{\ssc})$. Thus $\varpi_1,\dots,\varpi_l
\in X^*(A_{\ssc})$ are the fundamental weights for our root system. 
We now extend each fundamental weight $\varpi_i$ to a character $\omega_i$ on
$A$. (This involves a choice, since such an extension is well-defined only
modulo elements in $X^*(D)$.)   
 Additionally, we choose a basis $\omega_{l+1},\dots,\omega_n$ for
$X^*(D)$, where $n$ denotes the dimension of $A$. The characters
$\omega_1,\dots,\omega_n$ then form a $\mathbb Z$-basis for $X^*(A)$. 

As an example, consider the case of $GL_n$ with the standard choice for
$B=AU$. We identify both $X_*(A)$ and $X^*(A)$ with $\mathbb Z^n$. 
We also order the simple roots $\alpha_1,\dots,\alpha_{n-1}$ in the usual
way, so that
$\langle
\alpha_i,\nu \rangle=\nu_i-\nu_{i+1}$ for all $\nu \in X_*(A)$.  Then for
$i=1,\dots, n$  we may take
$\omega_i$ to be the weight
$1^i0^{n-i}
\in
\mathbb Z^n=X^*(A)$, the notation $1^i0^{n-i}$ meaning that $1$ is repeated
$i$ times and that $0$ is repeated $n-i$ times. Notice that for $\nu \in
X_*(A)$ the  expression $\nu_1+\dots +\nu_i$
occurring earlier is equal to $\langle \omega_i,\nu \rangle$. 

We will need the standard parabolic subgroups $P=MN$ of $G$, standard
meaning that $P$ contains $B$. Here $N$ is the unipotent radical of $P$, and
we always assume that the Levi subgroup
$M$ is the unique one containing $A$. The set $\Delta$ decomposes as the
disjoint union of $\Delta_M$ and $\Delta_N$, where $\Delta_M$ is the set of
simple roots of $M$ and $\Delta_N$ is the set of simple roots of $G$
occurring in the Lie algebra of $N$. We write
$A_M$ for the identity component of the center of $M$. Thus $A_M$ is a
subtorus of
$A$. 

We will need the real vector space $\mathfrak a:=X_*(A)\otimes
_\mathbb Z \mathbb R$. Sometimes we
will also need its $\mathbb Q$-subspace $\mathfrak a_\mathbb Q
:=X_*(A)\otimes _\mathbb Z
\mathbb Q$. We say that $x \in \mathfrak a$ is \emph{dominant} if $\langle
\alpha,x
\rangle
\ge 0$ for all
$\alpha
\in
\Delta$, and we write $\mathfrak a_{\dom}$ for the set of dominant elements
in $\mathfrak a$. We also need the usual partial order on $\mathfrak a$:
for $x,y \in \mathfrak a$ we say  that $x \le y$ if and only if $y-x$ is a
non-negative linear combination of simple coroots. In the case of
$GL_n$ an element
$\nu
\in
\mathfrak a=\mathbb R^n$ is dominant if and only if $\nu_1 \ge \dots \ge
\nu_n$, and the partial order is the same as the one discussed earlier (but
for $\mathbb R^n$ instead of $\mathbb Q^n$).  
Later we will also need the analogous partial order on the dual space
$\mathfrak a^*$: for $x,y  \in \mathfrak a^*$ we say that $x \le y$ if and
only if $y-x$ is a non-negative linear combination of simple roots. 

We write $W$ for the Weyl group of $A$ in $G$ and choose a $W$-invariant
Euclidean inner product $(x,y)$ on $\mathfrak a$, which we occasionally use
to identify $\mathfrak a$ with its dual.

More generally we also need $\mathfrak
a_M:=X_*(A_M)\otimes_\mathbb Z \mathbb R$. The obvious inclusion $X_*(A_M)
\subset X_*(A)$ lets us view $\mathfrak a_M$ as a subspace of
$\mathfrak a$. In fact $\mathfrak a_M=\{x \in \mathfrak a: \langle
\alpha,x \rangle =0 \quad \forall \, \alpha \in \Delta_M \}$. The subspace
$\mathfrak a_M$ of $\mathfrak a$ has a natural complement, namely the span
of the simple coroots for $M$. Thus $\mathfrak a_M$ is in a natural way a
direct summand of $\mathfrak a$, and we write $p_M:\mathfrak a
\twoheadrightarrow
\mathfrak a_M$ for the projection map. Write $W_M$ for the Weyl group of $A$
in $M$. Identifying $W_M$ with a subgroup of $W$, we can also describe
$\mathfrak a_M$ as the set of fixed points of $W_M$ on $\mathfrak a$, and
$p_M$ as the map sending $x$ to the average $|W_M|^{-1}\sum_{w \in W_M}
wx$ of the points in the $W_M$-orbit of $x$. Using our inner product, we can
also view $p_M$ as orthogonal projection on $\mathfrak a_M$.

We write $\Lambda_M$ for the quotient of $X_*(A)$ by the coroot lattice for
$M$. Since we have assumed that the derived group of $G$ is simply
connected, the same is true for $M$, and therefore  $\Lambda_M$ is a
free abelian group, which we identify with the image of $X_*(A)$ under
$p_M$, so that $\Lambda_M$ becomes a lattice in $\mathfrak a_M$. 

Returning to the standard parabolic subgroup $P=MN$, we obtain an open cone
$\mathfrak a_P^+$ in $\mathfrak a_M$ by putting 
\[
\mathfrak a_P^+:=\{x \in \mathfrak a_M: \langle \alpha,x\rangle > 0 \quad \,
\forall \, \alpha \in \Delta_N \}.
\]
We then define $\Lambda_P^+$ to be the intersection (in $\mathfrak a_M$) of
the lattice $\Lambda_M$ and the cone $\mathfrak a_P^+$. Recall that
$\mathfrak a_{\dom}$ is the disjoint union of the cones $\mathfrak a_P^+$ as
$P$ varies through all standard parabolic subgroups. Finally we define a
discrete subset $\mathcal N_G$ of $\mathfrak a_{\dom}$ by 
\[
\mathcal N_G:=\coprod_P \Lambda_P^+
\]
where $P$ runs over all standard parabolic subgroups. 

Let us work out the example of $GL_n$. To get a standard parabolic subgroup
$P=MN$   we need a partition $n=n_1+\dots+n_r$ of $n$. 
Given $\nu \in \mathfrak a=\mathbb R^n$, we refer to its first $n_1$ entries
as the first batch, its next $n_2$ entries as the second batch, and so on.
There are $r$ batches in all.  Then
$\mathfrak a_M$ is the subspace of $\mathbb R^n$ consisting of $n$-tuples
such that all entries in the first batch are equal to each other, all
entries in the second batch are equal to each other, and so on. The lattice
$\Lambda_M$ consists of such vectors which satisfy the additional condition
that the common value $\bar\nu_j$ of the  entries in the $j$-th batch
satisfies
$n_j\bar \nu_j \in \mathbb Z$. Finally, $\Lambda_P^+$ consists of vectors
satisfying all these conditions as well as the additional condition that
$\bar\nu_1 >\bar\nu_2 > \dots > \bar\nu_r$. The reader should observe that
$\mathcal N_G$ is precisely the set of vectors arising as the slopes of the Newton
polygon of some monic polynomial of degree $n$.

\subsection{The closest point map $r:\mathfrak a \to \mathfrak a_{\dom}$}
With these definitions out of the way, we now turn to a general root
theoretic construction which produces Newton polygons in the special case of
$GL_n$. In fact this construction is a standard one, one viewpoint being
that it is the most elementary case of Langlands' combinatorial lemma,
another being that it is the root theoretic fact needed for the Langlands
classification. An excellent textbook reference is \cite{knapp86}. The
basic idea is as follows. Let $x \in \mathfrak a$. Then there is a unique
point  $y
\in \mathfrak a_{\dom}$ that is closest to $x$ (in the metric on $\mathfrak
a$ obtained from our Euclidean inner product). In the case of $GL_n$ the map
$x \mapsto y$ is the map needed to produce the Newton polygon from the data
of the valuations of the coefficients of the polynomial $f(T)$. (Strictly
speaking, this is only true when every coefficient is non-zero.) 

The next
proposition summarizes all we need to know about $x \mapsto y$.
Much of what we need is in \cite{knapp86}; the rest will be verified later
(see
\ref{sub.pf131}).
 
\begin{proposition}\label{prop.131}
  \hfill 
\begin{enumerate}
\item  For all $x \in \mathfrak a$ there exists a unique point
$r(x)
\in
\mathfrak a_{\dom}$ that is closest to $x$.  The map $r$ is a continuous
retraction of $\mathfrak a$ onto
$\mathfrak a_{\dom}$. 
\item Let $y \in \mathfrak a_{\dom}$, and let $P=MN$ be the unique standard
parabolic subgroup such that $y \in \mathfrak a_P^+$. Then
\[
r^{-1}(y)=\{ x \in \mathfrak a:\text{$p_M(x)=y$ and $x \le y$} \}.
\] 
\item The point $r(x)$ can also be characterized as the unique minimal 
element in the set $\{\mu \in \mathfrak a_{\dom}:x \le \mu \}$. In other
words, $x \le r(x)$, and if $\mu \in \mathfrak a_{\dom}$ and $x
\le
\mu$, then
$r(x)
\le
\mu$. 

\item $r(\mathfrak a_\mathbb Q)=\mathfrak a_\mathbb Q \cap \mathfrak
a_{\dom}$.
\item $r(X_*(A))=\mathcal N_G$.
\end{enumerate}
\end{proposition}

 We actually need a slight generalization, since we
will sometimes be taking the valuation of $0$, which will yield $-\infty$.
To handle this it is best to proceed a bit differently. We will now use
the map $x \mapsto (\langle \omega_1,x\rangle,\dots,\langle
\omega_n,x\rangle)$ to identify $\mathfrak a$ with $\mathbb R^n$. Put
$\tilde{\mathbb R}:=\mathbb R \cup \{-\infty\}$ with the obvious topology,
for which a neighborhood base at $-\infty$ is provided by the intervals
$[-\infty,a)$. Then, as we will check later (see
\ref{sub.cont}) the map $r$ extends continuously to a map, still called $r$,
from 
$\tilde{\mathbb R}^l\times \mathbb R^{n-l}$ to $\mathfrak a_{\dom}$. 

\begin{proposition}\label{prop.132} \hfill
\begin{enumerate}
\item The image under $r$ of $\tilde{\mathbb Z}^l\times \mathbb Z^{n-l}$ 
is $\mathcal N_G$. 
\item Let $\mu \in \mathfrak a_{\dom}$ and $d \in \tilde{\mathbb R}^l\times
\mathbb R^{n-l}$. Then $r(d) \le \mu$ if and only if $d_i \le \langle
\omega_i,\mu \rangle$ for $i=1,\dots,l$ and $d_i = \langle
\omega_i,\mu \rangle$ for $i=l+1,\dots,n$. 
\item Let $\mu \in \mathfrak a_{\dom}$ and $d \in \tilde{\mathbb R}^l\times
\mathbb R^{n-l}$. Then $r(d) = \mu$ if and only if 
\begin{equation}
\begin{split}\label{eq.rev.vers} 
&d_i \le \langle \omega_i,\mu \rangle \text{ for all $i \in I_\mu$, and}\\
&d_i = \langle \omega_i,\mu \rangle \text{ for all $i \notin I_\mu$},
\end{split}
\end{equation}
\end{enumerate} 
where $I_\mu:=\{i \in \{1,\dots,l\}: \langle \alpha_i,\mu \rangle =0
\}$.
\end{proposition} 

This will be proved in subsection \ref{sub.pf132}. 

\subsection{Newton strata for the adjoint quotient $\mathbb A=A/W$}
We will be interested in the quotient variety $A/W$. The ring $R$ of regular
functions on $A/W$ is the
$F$-algebra obtained as the subring of $W$-invariant elements in the group
algebra $F[X^*(A)]$ of $X^*(A)$. When we view $\lambda \in X^*(A)$ as an
element of $F[X^*(A)]$ we denote it by $e^\lambda$. For
$i=1,\dots,n$ we put
$c_i:=\sum_{\lambda \in W\omega_i} e^\lambda$ (with $W\omega_i$ denoting
the Weyl group orbit of $\omega_i$ in $X^*(A)$). For $i=l+1,\dots,n$ we
then have simply $c_i=e^{\omega_i}$. It is well-known (see
\cite{bourbaki.root}) that the monomials
$c_1^{d_1}\dots c_n^{d_n}$ with $d_i \in \mathbb Z$ for all $i$ and $d_i \ge
0$ for $i=1,\dots,l$ form a basis for the $F$-vector space $R$. In other
words $R$ is the ring $F[c_1,\dots,c_n,c_{l+1}^{-1},\dots,c_n^{-1}]$
obtained from the polynomial ring $F[c_1,\dots,c_n]$ by inverting
$c_{l+1},\dots,c_n$, and the morphism 
\[
c:A \to \mathbb A^l \times \mathbb
G_m^{n-l}
\] 
defined by $c(a)=(c_1(a),\dots,c_n(a))$ identifies
$A/W$ with
$\mathbb A^l \times \mathbb G_m^{n-l}$. From now on we will write $\mathbb
A$ for $A/W$, even though it is only an affine space when $G=G_{\ssc}$. 

As an example, consider the case of $GL_n$, again taking $\omega_i$ to
be the weight
$1^i0^{n-i}
\in
\mathbb Z^n=X^*(A)$.  Writing $(a_1,\dots,a_n)$
for the $n$ diagonal entries of $a \in A(\bar F)$, we then have
\[
c_i(a)=\sum_I \prod_{j \in I} a_j,
\]
in which the sum is taken over all subsets $I$ of $\{1,\dots,n\}$ having
cardinality $i$. In other words $c_i$ is the usual $i$-th elementary
symmetric function of $a_1,\dots,a_n$, so that the characteristic
polynomial of $a$ is $T^n-c_1T^{n-1}+ \dots +(-1)^n c_n$. On the other hand
the eigenvalues  $a_1,\dots,a_n$ of $a$  are
obtained by evaluating certain characters of $A$ on $a$.  

We now return to the general case. Let $a \in A(\bar F)$. 
The example of $GL_n$ suggests 
that our goal should be to relate the valuations of the
elements 
$c_i(a)$ to those of 
$\lambda(a)$ (for characters
$\lambda \in X^*(A)$). 
To keep track of the numbers $\val\lambda(a)$ we
define an element 
$\nu_a \in \mathfrak a_\mathbb Q$ by requiring that
\[
\langle \lambda,\nu_a \rangle =\val \lambda(a) \qquad \forall \, \lambda \in
X^*(A).
 \] 
Thus $a \mapsto \nu_a$ is a
surjective homomorphism from $A(\bar F)$ to $\mathfrak a_\mathbb Q$. 
Clearly $\nu_{wa}=w\nu_a$ for any $w \in W$. 

Let $c = (c_1,\dots,c_n) \in \bar F^l \times (\bar F^\times)^{n-l}=\mathbb
A(\bar F)$. Put $\tilde {\mathbb Q}:=\mathbb Q \cup \{-\infty\}$ and define
an $n$-tuple $d_c \in \tilde {\mathbb Q}^l \times \mathbb Q^{n-l}$ by 
\[
d_c:=(\val c_1,\dots,\val c_n).
\]
The next result, which will be proved in section \ref{sec.pf141}, shows how
to recover the
$W$-orbit of
$\nu_a$ from
$d_{c(a)}$. 

\begin{theorem}\label{thm.rnu}
Let $a \in A(\bar F)$. Then $r(d_{c(a)})$ is the unique dominant element in
the $W$-orbit of $\nu_a$. 
\end{theorem} 

We are more interested in $\mathbb A(F)$ than in $\mathbb A(\bar F)$, and
for $\mu \in \mathcal N_G$ we now put 
\begin{align}
\mathbb A(F)_\mu:&=\{ c \in \mathbb A(F) : r(d_c)=\mu \} \\
\mathbb A(F)_{ \le \mu}:&=\{ c \in \mathbb A(F) : r(d_c) \le \mu \}.
\end{align} 
The sets $\mathbb A(F)_\mu$ are the \emph{Newton strata} referred to in the
title of this article. 

\begin{theorem}\label{thm.mu.strat}\hfill 
\begin{enumerate}
\item The sets $\mathbb A(F)_\mu$ are non-empty and 
\[
\mathbb A(F)=\coprod_{\mu \in \mathcal N_G} \mathbb A(F)_\mu.
\]
\item $\mathbb A(F)_{\le \mu}=\coprod_{\{\nu \in \mathcal N_G:\nu \le \mu\}}
\mathbb A(F)_\nu$.
\item $\mathbb A(F)_{\le \mu}$ consists of $n$-tuples $(c_1,\dots,c_n) \in
F^l \times (F^\times)^{n-l}=\mathbb A(F)$ such that 
\begin{align}
\val c_i &\le \langle \omega_i,\mu \rangle \quad \text{ for $i=1,\dots,l$,
and }  \\
\val c_i &= \langle \omega_i,\mu \rangle \quad \text{ for $i=l+1,\dots,n$.}
\end{align} 
\item $\mathbb A(F)_{ \mu}$ consists of $n$-tuples $(c_1,\dots,c_n) \in
F^l \times (F^\times)^{n-l}=\mathbb A(F)$ such that 
\begin{align}
\val c_i &\le \langle \omega_i,\mu \rangle \quad \text{ for all $i\in
I_\mu$, and }  \\
\val c_i &= \langle \omega_i,\mu \rangle \quad \text{ for
all $i\notin I_\mu$,}
\end{align} 
where $I_\mu:=\{i \in \{1,\dots,l\}: \langle \alpha_i,\mu \rangle =0
\}$.
\item Let $c \in \mathbb A(F)$ and choose $a \in A(\bar F)$ such that
$c(a)=c$. Then $c \in \mathbb A(F)_\mu$ if and only if $\mu$ is the unique
dominant element in the $W$-orbit of $\nu_a$. 
\end{enumerate}
\end{theorem}
\begin{proof}
(1) follows from the first part of Proposition \ref{prop.132}, and (2) is
clear from the definitions. (3) and (4) follow from the second and third
parts of Proposition
\ref{prop.132}. Finally, (5) follows from Theorem \ref{thm.rnu}. 
\end{proof}

\subsection{Dimensions of Newton strata}
Let us agree to assign dimensions (over $k$) to fractional ideals
$P^j:=\pi^j\mathcal O$ ($j \in \mathbb Z$) as follows: 
\[
\dim P^j = -j.
\] 
This definition is reasonable since we then have
\[
\dim P^i - \dim P^j=j-i,
\]
and, when $j \ge i$, we can interpret $j-i$ as the $k$-dimension of
$P^i/P^j$ if we use Greenberg's method \cite{greenberg61} to view $P^i/P^j$
as the set of $k$-points of a variety over $k$. 

Similarly, we assign to the subset $P^{j_1}\times\dots \times P^{j_l}$ of
$F^l$ the dimension $-(j_1 + \dots + j_l)$. Finally, for any coset $C$ of
$(\mathcal O^\times)^{n-l}$ in $(F^\times)^{n-l}$, we assign to the subset
$P^{j_1}\times\dots \times P^{j_l} \times C$ of $F^l
\times (F^\times)^{n-l}=\mathbb A(F)$ the dimension $-(j_1 + \dots + j_l)$.
This is compatible with the heuristic idea that the codimension of 
$P^{j_1}\times\dots \times P^{j_l} \times C$ in $P^{i_1}\times\dots \times
P^{i_l} \times C$ should be $(j_1+\dots+j_l)-(i_1+\dots +i_l)$ when $i_\beta
\le j_\beta $ for $\beta=1,\dots,l$. 

We see from the third part of Theorem \ref{thm.mu.strat} that the subset
$\mathbb A(F)_{\le \mu}$ of $\mathbb A(F)=F^l \times (F^\times)^{n-l}$ 
is equal to 
\[
P^{-[\langle \omega_1,\mu \rangle]} \times \dots \times P^{-[\langle
\omega_l,\mu \rangle]} \times C
\]
 for a suitable coset $C$ of $(\mathcal O^\times)^{n-l}$ in
$(F^\times)^{n-l}$.  Here $[a]$ denotes the greatest integer in $a$ (for $a
\in \mathbb R$).  Therefore the dimension of $\mathbb A(F)_{\le \mu}$ is
given by 
\begin{equation}\label{eq.dim1}
\dim \mathbb A(F)_{\le \mu} = \sum_{i=1}^l [\langle \omega_i,\mu \rangle].
\end{equation}
These numbers are not particularly meaningful since they depend on the
extensions $\omega_i$ of $\varpi_i$ we chose in \ref{subsec.rtp}. However,
given $\mu,\nu \in \mathcal N_G$ with $\nu \le \mu$, the set $\mathbb
A(F)_{\le\nu}$ is contained in $\mathbb A(F)_{\le\mu}$ with $k$-codimension 
\begin{equation}\label{eq.codim1}
\sum_{i=1}^l [\langle \omega_i,\mu \rangle]-\sum_{i=1}^l [\langle
\omega_i,\nu \rangle],
\end{equation}
a non-negative integer independent of the choice of extensions $\omega_i$. 
(To verify the independence statement use the fact that for any $\lambda
\in X^*(D)$ and
$\nu
\in
\mathcal N_G$ the number
$\langle\lambda,\nu\rangle$ is an integer.)

\subsection{Comparison with Chai's work}
An expression like $\sum_{i=1}^l [\langle \omega_i,\mu \rangle]$ appears in
Chai's \cite{chai00} conjectural formula for the dimensions of Newton strata
in Shimura varieties and again in Rapoport's conjectural formula
\cite{rapoport05} for the dimensions of affine Deligne-Lusztig varieties
inside  the affine Grassmannian. It is intriguing that it arises in such a
simple way in the context of the Newton stratification on $\mathbb A(F)$. 

Before comparing the codimensions \eqref{eq.codim1} with those in Chai's
article \cite{chai00} let us rewrite \eqref{eq.codim1} slightly. We have
extended the fundamental weights $\varpi_i$ to characters $\omega_i \in
X^*(A)$. There is another way to extend the fundamental weights to
``characters'' on $A$, more precisely, to elements in $X^*(A)_\mathbb Q$.
For this we use the direct sum decomposition
$X^*(A)_{\mathbb Q}=X^*(A_{\ssc})_\mathbb Q \oplus X^*(A_G)_\mathbb Q$ to
view the fundamental weights $\varpi_i$ as elements of $X^*(A)_\mathbb Q$
that happen to lie in the lattice $X^*(A_{\ssc})$ in the direct summand
$X^*(A_{\ssc})_\mathbb Q$ of $X^*(A)_\mathbb Q$. In other words, for
$i=1,\dots,l$ we may view $\varpi_i$ as the unique element of
$X^*(A)_\mathbb Q$ satisfying the conditions $\langle
\varpi_i,\alpha^\vee_j \rangle =\delta_{ij}$ for $1 \le j \le l$ and
$\langle \varpi_i,\nu \rangle=0$ for all $\nu \in \mathfrak a_G$. 

Now consider a dominant coweight $\mu \in X_*(A)$ and an element $\nu \in
\mathcal N_G$ such that $\nu \le \mu$, so that  in particular $\mu-\nu \in
X_*(A_{\ssc})_\mathbb Q$. We are going to rewrite \eqref{eq.codim1}. Since
$\mu \in X_*(A)$, the numbers $\langle \omega_i,\mu \rangle$ are integers,
from which it follows that
\[
[\langle \omega_i,\mu \rangle]-[\langle \omega_i,\nu \rangle]=-[\langle
\omega_i,\nu-\mu \rangle]=\lceil \langle \omega_i,\mu-\nu \rangle \rceil
\]
where $\lceil r \rceil$ denotes the smallest integer that is greater than
or equal to the real number $r$. Moreover, since $\mu-\nu \in
X_*(A_{\ssc})_\mathbb Q$, we have $\langle \omega_i,\mu-\nu \rangle=\langle
\varpi_i,\mu-\nu \rangle$. We conclude that for dominant  $\mu \in X_*(A)$
and  $\nu \in \mathcal N_G$ such that $\nu \le \mu$, the set $\mathbb
A(F)_{\le \nu}$ is contained in $\mathbb A(F)_{\le \mu}$ with
$k$-codimension 
\begin{equation}\label{eq.codim2}
\sum_{i=1}^l \lceil \langle \varpi_i,\mu-\nu \rangle \rceil.
\end{equation}
Comparing this with the last displayed formula in \cite[Thm.~7.4]{chai00},
one sees that the codimension of $\mathbb
A(F)_{\le \nu}$  in $\mathbb A(F)_{\le \mu}$ coincides with the value
predicted by Chai for the codimension of the Newton stratum indexed by
$\nu$ in any Shimura variety with associated minuscule coweight $\mu$ (with
$G$ split over $\mathbb Q_p$). 

Rapoport's conjectural formula (see \cite{rapoport05},
\cite[eqn.~(0.2)]{mierendorff05.pre}) for dimensions of affine
Deligne-Lusztig varieties inside  affine Grassmannians can be given a
similar treatment. This has been done implicitly in \cite{ghkr.pre}, using
the notion of \emph{defect} that we will be discussing in the remainder of
this introduction.

\subsection{Definition of $d_G:\mathcal N_G \to \mathbb Q$}
For a rational number $r$ let us write $\fr(r) \in [0,1)$ for the
fractional part of $r$; thus
\[
[r]=r-\fr(r).
\]
One sees immediately that
\[
\sum_{i=1}^l [\langle \omega_i,\nu \rangle]=\langle \rho',\nu
\rangle-\sum_{i=1}^l \fr \langle \omega_i,\nu \rangle,
\]
where $\rho':= \sum_{i=1}^l \omega_i$, an extension to $A$ of the usual
character 
\[
\rho:=\sum_{i=1}^l\varpi_i=\frac{1}{2}\sum_{\alpha >0} \alpha \in
X^*(A_{\ssc}).
\]

We remind the reader that for any $\lambda \in X^*(D)$ and $\nu \in \mathcal
N_G$ the number
$\langle\lambda,\nu\rangle$ is an integer. Therefore
$\fr\langle\omega_i,\nu\rangle$ is independent of the choice of extension
$\omega_i$ of $\varpi_i$, and, moreover, $\fr\langle \omega_i,\nu
\rangle=0$ for $i=l+1,\dots,n$. Thus we may define a function 
\[
d_G:\mathcal N_G \to  \mathbb Q,
\]
independent of the choice of extensions $\omega_i$, by putting
\[
d_G(\nu):=\sum_{i=1}^l \fr \langle \omega_i,\nu \rangle=\sum_{i=1}^n \fr
\langle \omega_i,\nu \rangle.
\]

Using this definition, we may rewrite our formula \eqref{eq.dim1} as
\begin{equation}
\dim \mathbb A(F)_{\le \nu} = \langle \rho',\nu \rangle - d_G(\nu).
\end{equation}

\subsection{Computation of $d_G(\nu)$}  We are going to calculate
$d_G(\nu)$ for all $\nu \in \mathcal N_G$.  Since $\Lambda_G^+=\Lambda_G$, 
the set $\mathcal N_G$ contains the lattice
$\Lambda_G$, and it turns out that the main point is to
calculate $d_G(\nu)$ when $\nu \in \Lambda_G$. Indeed the general case may
be reduced to this special one by means of the next lemma.

\begin{lemma}\label{lem.dgm} 
Let $\nu \in \mathcal N_G$, and let $P=MN$ be the unique standard parabolic
subgroup such that $\nu \in \Lambda_P^+$. Then there is an
equality
\[
d_G(\nu)=d_M(\nu),
\]
where on the right side $\nu$ is  viewed as an element of $\Lambda_M
\subset \mathcal N_M$ and $d_M:\mathcal N_M \to \mathbb Q$ is the analog
for $M$ of $d_G$. 
\end{lemma}
\begin{proof}  
 The derived group of $M$ is simply connected, and one sees immediately
that the  characters $\omega_i$ play the same role for $M$ as they do for
$G$. Therefore the $n$ numbers $\fr\langle \omega_i,\nu \rangle$ we get for
$M$ are the same as the $n$ numbers we get for $G$.
\end{proof}

We are now reduced  to computing $d_G(\nu)$ when $\nu \in \Lambda_G \subset
\mathcal N_G$.  
For this we need to define, for each such $\nu$,  a non-negative integer
$\defect_G(\nu)$ called the \emph{defect} of $\nu$. Recall that the
extended affine Weyl group $\tilde W$ is the semidirect product $\tilde
W=X_*(A)\rtimes W$, and therefore comes with a surjective homomorphism 
\[
p:\tilde W \twoheadrightarrow W.
\] 
Since $W$ acts trivially on the quotient $\Lambda_G$ of $X_*(A)$, there is
also an obvious surjective homomorphism 
\[
q:\tilde W \twoheadrightarrow \Lambda_G,
\]
which is trivial on $W$ and on $X_*(A)$ is the canonical surjection $X_*(A)
\twoheadrightarrow \Lambda_G$. The kernel of $q$ is the affine Weyl group.
Inside $\tilde W$ we have the subgroup consisting of all elements that
preserve the set of simple affine roots (equivalently, that preserve the
base alcove in $\mathfrak a$), and $q$ restricts to an isomorphism from
this subgroup onto $\Lambda_G$; inverting this isomorphism we obtain the
usual section $s$ of $q$, a homomorphism 
\[
s:\Lambda_G \to \tilde W
\]
such that $qs=\id$. The composition $ps$ is a homomorphism 
\[
ps:\Lambda_G \to W.
\]
The defect $\defect_G(\nu)$ is then defined to be 
\begin{equation}
\defect_G(\nu):=\dim \mathfrak a-\dim \mathfrak a^{w_\nu},
\end{equation}
where $w_\nu:=ps(\nu)$. Here, for $w \in W$ we are using $\mathfrak a^w$
to denote the set of fixed points of $w$ in $\mathfrak a$. It may be worth
noting that the difference $\dim \mathfrak a-\dim \mathfrak a^w$ also
occurs in the dimension formula for affine Springer fibers (see
\cite{kazhdan-lusztig88} and \cite{bezrukavnikov96}).

\begin{theorem}\label{thm.d=def}
For all $\nu \in \Lambda_G \subset \mathcal N_G$ there is an equality 
\begin{equation}
d_G(\nu)=\frac{1}{2}\defect_G(\nu).
\end{equation}
In particular $d_G(\nu) \in \frac{1}{2} \mathbb Z$.
\end{theorem}
\begin{proof}
This will be proved in subsection \ref{subsec.4.2}.
\end{proof}

\subsection{$B(G)$ and $\defect_G(b)$} 
The set $\mathcal N_G$ also arises in the closely related context of
isocrystals with $G$-structure. As in
\cite{kottwitz85,rapoport-zink96,kottwitz97}, given a split group
$G$ over $\mathbb Q_p$ and a basic element $b \in G(L)$ (where $L$ is the
completion of a maximal unramified extension of $\mathbb Q_p$), we denote
by $J_b$ the inner form of $G$ obtained by twisting the action of the
Frobenius automorphism $\sigma$ of $L/\mathbb Q_p$ on $G(L)$ by $b$, so
that 
\begin{equation}\label{eq.J}
J_b(\mathbb Q_p)=\{g \in G(L) : g^{-1}b\sigma(g)=b\}.
\end{equation}
For any connected reductive group $H$ over $\mathbb Q_p$ we denote by
$\rk_{\mathbb Q_p}(H)$ the dimension of any maximal $\mathbb Q_p$-split
torus of $H$. Recall from \cite{kottwitz85} that there is a bijection $b
\mapsto \nu_b$ from the set of $\sigma$-conjugacy classes of basic elements
 $b \in G(L)$ to the abelian group $\Lambda_G$. 

\begin{theorem}\label{thm.bdef}
Let $\nu \in \Lambda_G \subset \mathcal N_G$ and let $b$ be any basic
element in $G(L)$ such that $\nu_b=\nu$. Then $\rk_{\mathbb Q_p}(G)=\dim
\mathfrak a$ and $\rk_{\mathbb Q_p}(J_b)=\dim \mathfrak a^{w_\nu}$, so
that 
\begin{equation}
\defect_G(\nu)=\rk_{\mathbb Q_p}(G)-\rk_{\mathbb Q_p}(J_b).
\end{equation}
\end{theorem}
\begin{proof}
This will be proved in subsection \ref{subsec.4.3}.
\end{proof}

More generally (see \cite{kottwitz85,rapoport-richartz96}), provided that
the derived group of $G$ is simply connected, as we once again assume, there
is a canonical map $b \mapsto \nu_b \in \mathcal N_G$ which induces a
bijection from $B(G)$ (the set of $\sigma$-conjugacy classes in $G(L)$) to
the set $\mathcal N_G$. For $b \in G(L)$ one can still define a group $J_b$
(again see \cite{kottwitz85,rapoport-zink96,kottwitz97}) satisfying
\eqref{eq.J}. If $\nu_b \in \Lambda_P^+ \subset \mathcal N_G$ for standard
$P=MN$, then  $b$ is $\sigma$-conjugate to
a basic element $b_M \in M(L)$ such that $\nu_{b_M}=\nu_b$. Moreover, $J_b$
is isomorphic to $J_{b_M}$,  an inner form of $M$.  Combining Lemma
\ref{lem.dgm} and Theorems \ref{thm.d=def} and \ref{thm.bdef}, and noting
that $G$ has the same $\mathbb Q_p$-rank as $M$, we
obtain the 
\begin{corollary}
For any $b \in G(L)$ there is an equality
\begin{equation}
d_G(\nu_b)=\frac{1}{2}\defect_G(b),
\end{equation}
where $\defect_G(b)$ is defined by 
\begin{equation}
\defect_G(b):=\rk_{\mathbb Q_p}G-\rk_{\mathbb Q_p}J_b.
\end{equation}
\end{corollary}
This corollary was used implicitly in \cite{ghkr.pre} in order to
reformulate the dimension conjecture of Rapoport in terms of
$\defect_G(b)$. 

\subsection{What if the derived group of $G$ is not simply connected?} 
The best way to proceed in this case is probably to choose a
central extension $G'$ of $G$ by a split torus, doing this in such a way
that the derived group of $G'$ is simply connected. For example one would
study
$PGL_n$ by relating it to $GL_n$. 

\section{Proofs of Propositions \ref{prop.131} and \ref{prop.132}}
\subsection{Proof of Proposition \ref{prop.131}}\label{sub.pf131}
(1) For any non-empty closed convex subset $C$ in the Euclidean space
$\mathfrak a$ and any $x \in \mathfrak a$, there exists a unique point $r(x)
\in C$ closest to $x$. Obviously $r:\mathfrak a \to C$ is a retraction, and
$r$ is continuous because it satisfies the inequality
\[
d(r(x_1),r(x_2)) \le d(x_1,x_2),
\]
where $d(-,-)$ denotes the Euclidean metric on $\mathfrak a$. To prove this
inequality we may replace $C$ by the line segment joining $r(x_1)$ and
$r(x_2)$, and then the proof is easy.

(2) This is simply a reformulation of Lemma 8.56 in \cite{knapp86}. 

(3) Uniqueness is clear since the partial order $\le$ has the property
that $x \le x'$ and $x' \le x$ imply $x=x'$. We see from (2) that $x \le
r(x)$. It remains to prove that if $\mu \in \mathfrak a_{\dom}$ and $x \le
\mu$, then $r(x) \le \mu$. Let $P=MN$ be the unique standard parabolic
subgroup such that $r(x) \in \mathfrak a_P^+$. We claim that
\[
r(x) = p_M(x) \le p_M(\mu) \le \mu.
\]
Indeed, $r(x)=p_M(x)$ follows from (2). The inequality $p_M(x) \le
p_M(\mu)$ follows from the inequality $x \le \mu$ since $p_M$ preserves
$\le$. Finally, the inequality $p_M(\mu) \le \mu$ is a consequence of the
dominance of $\mu$ (since $p_M(\mu)$ is the average of the points $w\mu$ ($w
\in W_M$) and the dominance of $\mu$ implies that $w\mu \le \mu$). 

(4) The map $r$ carries $\mathfrak a_\mathbb Q$ into  $\mathfrak a_\mathbb
Q \cap \mathfrak a_{\dom}$ since it follows from (2) that there exists
standard $P=MN$ such that $r(x)=p_M(x)$. Since $r$ is a retraction, it maps 
$\mathfrak a_\mathbb Q$ onto  $\mathfrak a_\mathbb
Q \cap \mathfrak a_{\dom}$.

(5) Let $y \in \mathfrak a_{\dom}$ and suppose that $y \in \mathfrak a_P^+$
for  standard $P=MN$. Then, using (2), we see that $y \in r(X_*(A))$
$\iff$ $p_M^{-1}(y) \cap \{x \in \mathfrak a: x \le y \}$ meets $X_*(A)$.
This happens $\iff$ $p_M^{-1}(y)$ meets $X_*(A)$ $\iff$ $y \in \mathcal
N_G$. 

\subsection{Continuity of $r:\tilde {\mathbb R}^l\times \mathbb R^{n-l} \to
\mathfrak a_{\dom}$}\label{sub.cont} 
 The reader might find it helpful to draw a picture for $G=SL_3$, in which
case it is clear that $r:\mathbb R^n \to \mathfrak a_{\dom}$ extends
continuously to   
$\tilde {\mathbb R}^l\times \mathbb R^{n-l} \to
\mathfrak a_{\dom}$. 

Before handling the general case, we pause to transcribe part (2) of
Proposition \ref{prop.131} using our identification $\mathfrak a=\mathbb
R^n$. The statement becomes: for standard $P=MN$ and $y \in \mathfrak
a_P^+$, the fiber $r^{-1}(y)$ is the set of $d=(d_1,\dots,d_n) \in \mathbb
R^n$ such that 
\begin{equation}
\begin{split}\label{eq.domg} 
&d_i \le \langle \omega_i,y \rangle \text{ for all $i \in I_M$, and}\\
&d_i = \langle \omega_i,y \rangle \text{ for all $i \notin I_M$},
\end{split}
\end{equation}
where $I_M:=\{i \in \{1,\dots,l\}: \alpha_i \in \Delta_M \}$. 

Now we prove that $r$ extends continuously. The proof will use a map $d
\mapsto d'$ from $\mathbb R^n$ to itself that we will now define.  We begin
with the special case in which 
$G=G_{\ssc}$, so that $n=l$ and $\omega_i=\varpi_i$ ($1 \le i \le l$). It
is now convenient to use our $W$-invariant inner product to identify
$\mathfrak a^*$ with $\mathfrak a$, so that $\omega_i$ can be regarded as
an element of $\mathfrak a$.  Since any
$y
\in
\mathfrak a_P^+$ is a strictly positive linear combination of $\omega_i$
($i \notin I_M$), the numbers
$\langle
\omega_i,y
\rangle$ occurring in \eqref{eq.domg} are non-negative for all $i$ and
strictly positive when $i \notin I_M$, and therefore we see that $r(d)=y$
$\iff$
$r(d')=y$, where $d'=(d'_1,\dots,d'_l) \in \mathbb R^l$ is defined by 
\[
d_i':=\max\{d_i,0\}.
\]

In the general case we write $\mathfrak a=\mathbb R^n$ as $\mathfrak
a_{\ssc}
\oplus
\mathfrak a_G$, where $\mathfrak a_{\ssc}$ is the analog of $\mathfrak a$
for the group $G_{\ssc}$. For $d=(d_{\ssc},z) \in \mathfrak a_{\ssc} \oplus
\mathfrak a_G$ we define $d' \in \mathfrak a$ by $d':=(d'_{\ssc},z)$, where
$d_{\ssc} \mapsto d'_{\ssc}$ is the map we just defined from $\mathfrak
a_{\ssc}=\mathbb R^l$ to itself. Again we have $r(d)=y$ $\iff$
$r(d')=y$. In other words $r:\mathbb R^n \to \mathfrak a_{\dom}$ is the
composition $rt$, where $t$ is the map $d \mapsto d'$. It is clear that $t$
extends continuously to a map 
\[
\tilde{\mathbb R}^l\times \mathbb R^{n-l} \to \mathbb R^n.
\]
Indeed, for this we may as well assume that $G=G_{\ssc}$, in which case 
the desired extension of $t$ is still given by $d \mapsto d'$, with
$\max\{d_i,0\}$ interpreted as $0$ when $d_i=-\infty$. Therefore $r$
extends continuously.

\subsection{Proof of Proposition \ref{prop.132}} \label{sub.pf132} 
As a byproduct of the proof we just gave, we see that for our extended map
$r:\tilde {\mathbb R}^l\times \mathbb R^{n-l} \to
\mathfrak a_{\dom}$ and $y \in \mathfrak a_P^+$, the fiber $r^{-1}(y)$ is
given by the set of
$d \in \tilde {\mathbb R}^l\times \mathbb R^{n-l}$ satisfying
\eqref{eq.domg}. This proves part (3) of Proposition \ref{prop.132}, and it
remains only to prove parts (1) and (2). 

(1) We already know from Proposition \ref{prop.131} that $r(\mathbb
Z^n)=\mathcal N_G$, so we just need to see that allowing some coordinates
of $d \in \tilde {\mathbb Z}^l\times \mathbb Z^{n-l}$ to be $-\infty$ does
not give anything new. This is clear from the discussion in subsection
\ref{sub.cont}, since $r(d)$ does not change when each $-\infty$ occurring
as a coordinate is replaced by some sufficiently negative integer. 

(2) Put $y:=r(d)$. First we prove ($\Longrightarrow$). So assume that $y
\le \mu$, which means that $\langle \omega_i,y \rangle \le \langle
\omega_i,\mu \rangle $ for $i=1,\dots,l$ and $\langle
\omega_i,y\rangle=\langle \omega_i,\mu \rangle $ for $i = l+1,\dots,n$.
This, together with \eqref{eq.domg} (which we know is true since we have
already proved part (3) of this proposition), shows that
\begin{equation}
\begin{split}\label{eq.domgd} 
&d_i \le \langle \omega_i,\mu \rangle \text{ for all $i =1,\dots,l$, and}\\
&d_i = \langle \omega_i,\mu \rangle \text{ for all $i =l+1,\dots,n$},
\end{split}
\end{equation}
as desired. 

Next we prove  ($\Longleftarrow$). Assume that \eqref{eq.domgd} holds.
Consider the point $d' \in \mathbb R^n$ constructed from $d$ in subsection
\ref{sub.cont}. Then \eqref{eq.domgd} holds with $d$ replaced by $d'$, or
in other words, $d'$ satisfies $d' \le \mu$. From part (3) of Proposition
\ref{prop.131} we see that $r(d') \le \mu$. Since $r(d)=r(d')$, this
concludes the proof. 

\section{Proof of Theorem \ref{thm.rnu}}\label{sec.pf141} 
Replacing $a$ by a suitable $W$-conjugate, we may assume that $\nu_a$ is
dominant. We must then show that $r(d_{c(a)})=\nu_a$. Let $P=MN$ be the
unique standard parabolic subgroup such that $\nu_a \in \mathfrak a_P^+$.
Then, appealing to part (3) of Proposition \ref{prop.132}, we see that we
must prove that 
\begin{align}
\val c_i(a) &\le \langle \omega_i,\nu_a \rangle \text{ for all $i \in I_M$,
and}\label{eq.p1}\\
\val c_i(a) &= \langle \omega_i,\nu_a \rangle \text{ for all $i \notin
I_M$,}\label{eq.p2} 
\end{align}  
where $I_M:=\{i \in \{1,\dots,l\}: \alpha_i \in \Delta_M \}$. 

For $i = l+1,\dots,n$ we have $c_i(a)=\omega_i(a)$, and therefore   
\eqref{eq.p2} holds for such $i$. Now suppose that $i$ is in the range
$i=1,\dots,l$. Then 
\[
c_i(a)=\sum_{\lambda \in W\omega_i} \lambda(a),
\]
so that 
\[
\val c_i(a) \le \max\{ \langle w\omega_i,\nu_a \rangle: w \in W\}.
\]
Since $\omega_i$ is dominant, $w\omega_i \le \omega_i$, and this implies
($\nu_a$ being dominant) that
\[
\langle w\omega_i,\nu_a \rangle \le 
\langle \omega_i,\nu_a \rangle.
\]
This proves \eqref{eq.p1} and shows that in order to prove \eqref{eq.p2}, it
is enough to show that, if $i \notin I_M$, then 
\begin{equation}\label{eq.hrt}
 \langle w\omega_i,\nu_a \rangle < 
\langle \omega_i,\nu_a \rangle
\end{equation}
unless $w\omega_i=\omega_i$. 

We already noted that $w\omega_i\le \omega_i$, which means that 
 \begin{equation}\label{eq.flat}
\omega_i-w\omega_i=\sum_{j=1}^l n_j \alpha_j 
\end{equation} 
for some non-negative integers $n_1,\dots,n_l$. 
Therefore the difference between the two sides of \eqref{eq.hrt} is
\begin{equation}\label{eq.numb}
\sum_{j=1}^l n_j\langle \alpha_j,\nu_a \rangle. 
\end{equation}
Since $\langle \alpha_j,\nu_a \rangle \ge 0$, all terms in this sum are
non-negative, and thus it will suffice to show that the $i$-th term is
strictly positive. Since $i \notin I_M$, we do know that $\langle
\alpha_i,\nu_a \rangle>0$, and it remains only to prove that $n_i > 0$. 

Recall our $W$-invariant inner product $(\cdot,\cdot)$ on $\mathfrak a$,
which we now use to identify $\mathfrak a$ with its dual; each simple root
$\alpha_j$ is then identified with a strictly positive multiple of the
corresponding simple coroot $\alpha_j^\vee$.  Taking the inner product of
each side of
\eqref{eq.flat} with
$\omega_i$, we see that
$n_i$ is strictly positive if and only if
\[
(\omega_i-w\omega_i,\omega_i)>0.
\]
For any pair $v$, $v'$ of distinct non-zero vectors having the same
Euclidean norm, the inner product $(v-v',v)$ is strictly positive
(Cauchy-Schwarz inequality). Applying this observation to the vectors
$\omega_i$, $w\omega_i$ completes the proof. 

\section{Proof of Theorems \ref{thm.d=def} and \ref{thm.bdef}} 
\subsection{A lemma about the extended affine Weyl group $\tilde W$}
\label{sub.lemW}
Recall from before the surjective homomorphisms $p:\tilde W
\twoheadrightarrow W$, $q:\tilde W \twoheadrightarrow \Lambda_G$, and the
section $s:\Lambda_G \to \tilde W$ of $q$. Write $V$ for $\mathfrak
a^*\otimes_\mathbb R \mathbb C$. We are interested in the reflection
representation of $W$ on $V$, which we then regard as an $n$-dimensional
representation of $\Lambda_G$ using the homomorphism $ps:\Lambda_G \to W$. 
The composition 
\[
X_*(A_G) \hookrightarrow X_*(A) \twoheadrightarrow \Lambda_G
\]
identifies $X_*(A_G)$ with a subgroup of finite index in $\Lambda_G$, and
this subgroup acts trivially on $V$; therefore $V$ must decompose as the
direct sum of $n$ $1$-dimensional representations of $\Lambda_G$. 
We are going to describe the $n$ characters of $\Lambda_G$ so obtained. 

For each $i=1,\dots,n$ we define a character $\chi_i:\Lambda_G \to \mathbb
C^\times$ as the following composition:
\[
\Lambda_G \hookrightarrow \mathfrak a_G \xrightarrow{\langle \omega_i,\cdot
\rangle} \mathbb R \twoheadrightarrow \mathbb R/\mathbb Z
\xrightarrow{\exp(2\pi i\cdot)} \mathbb C^\times.
\]
Of course $\chi_i$ is trivial for $i=l+1,\dots,n$, since $\langle
\omega_i,\cdot \rangle$ then takes integral values on the lattice 
$\Lambda_G$ in $\mathfrak a_G$.   

\begin{lemma}\label{lem.refl} 
As a $\Lambda_G$-module $V$ is isomorphic to the direct sum $\chi_1 \oplus
\dots \oplus \chi_n$ of the $1$-dimensional representations $\chi_i$. 
\end{lemma}
\begin{proof}
We already remarked that $V$ is really a representation of the quotient
$\Lambda_G/X_*(A_G)$ of $\Lambda_G$. Moreover each character $\chi_i$ is
clearly trivial on the subgroup $X_*(A_G)$. Therefore our problem really
concerns the quotient $\Lambda_G/X_*(A_G)$, which may be identified with a
subgroup of $P(R^\vee)/Q(R^\vee)$ (Bourbaki's notation), where $R^\vee$
is the coroot system of $G$, $P(R^\vee)$ is the lattice of weights for
$R^\vee$, and $Q(R^\vee)=X_*(A_{\ssc})$ is the coroot lattice for $G$. 

How do we express $\chi_i$ ($1 \le i \le l$) in terms of
$P(R^\vee)/Q(R^\vee)$? We use $\exp(2\pi i\cdot)$ to identify
$\mathbb R/\mathbb Z$ with the unit circle in $\mathbb C^\times$, and
$\mathbb Q/\mathbb Z$ with the group of complex roots of unity.  For
$\nu \in \Lambda_G \subset \mathfrak a_G$ pick $\mu \in X_*(A)$ such that
$p_G(\mu)=\nu$. We then claim that $\chi_i(\nu)$ is equal to $-\langle
\varpi_i,\bar\mu \rangle \in \mathbb Q/\mathbb Z \hookrightarrow \mathbb
C^\times$, where $\bar\mu$ denotes the image of $\mu $ under 
\[
X_*(A) \twoheadrightarrow \Lambda_G \twoheadrightarrow \Lambda_G/X_*(A_G)
\hookrightarrow P(R^\vee)/Q(R^\vee).
\] 
Indeed, since $\langle \omega_i,\mu \rangle \in \mathbb Z$, the rational
numbers $\langle \omega_i,\nu \rangle$ and $\langle \omega_i,\nu-\mu
\rangle $ become equal in $\mathbb Q/\mathbb Z$. Moreover, since $\nu-\mu
\in X_*(A_{\ssc})_\mathbb Q$, we have 
\[
\langle \omega_i,\nu-\mu \rangle = \langle \varpi_i,\nu-\mu
\rangle=-\langle \varpi_i,\mu \rangle. 
\]

Thus our entire problem can be reformulated in terms of the adjoint group
$G_{\ad}$ of $G$. So for the rest of this proof we will consider an
adjoint group $G$ and its extended affine Weyl group $\tilde W$. Thus
$\Lambda_G$ is now the finite abelian group $\Lambda:=P(R^\vee)/Q(R^\vee)$
and
$\chi_i$ ($i=1,\dots,l$) has become the character
$\chi_i:\Lambda \to \mathbb Q/\mathbb Z \hookrightarrow \mathbb
C^\times$ defined by
\[
\chi_i(\bar\nu)=-\langle \varpi_i,\nu \rangle,
\]
where $\nu\in P(R^\vee)$ is any representative for $\bar\nu \in \Lambda$. 

When $G$ decomposes as a direct product, our problem decomposes
accordingly, and therefore we may as well assume that $G$ is a simple
adjoint group. At this point we  use the classification of root systems.
All the information we need can be obtained easily from the tables of root
systems in Bourbaki \cite{bourbaki.root}. This is obviously the case for the
numbers
$\langle
\varpi_i,\nu \rangle$, but what about the decomposition of the
$\Lambda$-module $V$ as a sum of characters? 

We can also consider the representation $V'$ of $\tilde W$ on the space of
affine linear (complex-valued) functions on $\mathfrak a$. The simple
affine roots give a $\mathbb C$-basis for $V'$, and $\Lambda$
permutes this basis. Thus the representation on $V'$ of
$\Lambda$ is the permutation representation obtained from its
action on the set of vertices of the affine Dynkin diagram of $G$; this
information too is contained in Bourbaki. We have the short exact sequence 
\[
0 \to \mathbb C \hookrightarrow V' \to V \to 0,
\]
with $\mathbb C \hookrightarrow V'$ given by the constant functions on
$\mathfrak a$. Therefore the assertion we need to check (for each root
system) is that the set of $l+1$ characters of $\Lambda$ we get from its
permutation representation on the vertices of the affine Dynkin diagram
consists of the trivial character together with the $l$ characters
$\chi_i$. 

For example in type $A_{n-1}$ the group $\Lambda$ is
cyclic of order $n$ and permutes the $n$ vertices of the affine  Dynkin
diagram cyclically, so that each character of $\Lambda$ appears once
in this
$n$-dimensional permutation representation. In Bourbaki's numbering of the
fundamental weights, $\chi_1$ generates the character group of $\Lambda$ and
$\chi_i$ ($i=1,\dots,n-1$) is the $i$-th power of
$\chi_1$. 

As another example, for $E_6$ the group $\Lambda$ is cyclic of order $3$,
and has three orbits on the set of vertices of the affine Dynkin diagram,
two orbits having three elements, and one having one element. Therefore
$V'$ decomposes as a sum of two copies of each of the non-trivial
characters of $\Lambda$ and three copies of the trivial character. Again
using Bourbaki's numbering system, we find that $\chi_2=\chi_4$ are
trivial, $\chi_1=\chi_5$ are one of the two non-trivial characters, and
$\chi_3=\chi_6$ are the other non-trivial character. So it works! One
checks just as easily that it works for each simple root system. 

We remark that  since the reflection representation is defined over
$\mathbb Q$, it was only a matter of convenience that we used $\exp(2\pi
i\cdot)$ to identify $\mathbb Q/\mathbb Z$ with the group of complex roots
of unity. Any other identification of the two groups would have worked just
as well. 
\end{proof}

\subsection{Proof of Theorem \ref{thm.d=def}}\label{subsec.4.2} 
We will prove Theorem \ref{thm.d=def} by combining  Lemma
\ref{lem.refl}   with an argument of the same kind that Bezrukavnikov used
to prove
\cite[Lemma 3]{bezrukavnikov96}. 
 The characters $\chi_i$
($i=1,\dots,n$) introduced in subsection \ref{sub.lemW} take values in
$\mathbb Q/\mathbb Z \hookrightarrow \mathbb C^\times$. Theorem
\ref{thm.d=def} involves the fractional parts of $\langle \omega_i,\nu
\rangle$. We identify $\mathbb Q/\mathbb Z$ with the set of rational numbers
in the interval $[0,1)$, so that $\chi_i$ now takes values in $[0,1)$. With
this convention, we must show that for $\nu \in \Lambda_G$
\begin{equation}\label{eq.natu}
2\sum_{i=1}^n \chi_i(\nu) = \dim \mathfrak a - \dim \mathfrak a^{w},
\end{equation} 
where $w=w_\nu=ps(\nu)$. 
Since the reflection representation $V$ is self-contragredient, the
multiplicity of $\chi_i$ in $V$ is the same as that of $\chi_i^{-1}$.
Therefore the left side of \eqref{eq.natu} is equal to 
\[
\sum_{i=1}^n \chi_i(\nu)+\chi_i^{-1}(\nu). 
\]
Now note that $\chi_i(\nu)+\chi_i^{-1}(\nu)$ is $1$ if $\chi_i(\nu)\ne 0
\in \mathbb Q/\mathbb Z$ and is $0$ otherwise. Therefore, using Lemma
\ref{lem.refl}, we see that the left side of
\eqref{eq.natu} is $\dim \mathfrak a$ minus the number of times $1$ occurs
as an eigenvalue of $w$ acting on $\mathfrak a$, namely $\dim\mathfrak
a^w$. 

\subsection{Proof of Theorem \ref{thm.bdef}}\label{subsec.4.3} 
In this theorem, we are given $\nu \in \Lambda_G$, and we are
interested in basic elements $b \in G(L)$ such that $\nu_b=\nu$. There is a
unique $\sigma$-conjugacy class of such elements, and since the isomorphism
class of $J_b$ depends only on the $\sigma$-conjugacy class of $b$, it
suffices to prove the theorem for any particular $b$ we like. The most
convenient choice for $b$ is obtained as follows. The extended affine Weyl
group
$\tilde W$ can be identified with $N_G(A)(\mathbb Q_p)/A(\mathbb Z_p)$,
where $N_G(A)$ denotes the normalizer of $A$ in $G$. We take as $b$ any
element of $N_G(A)(\mathbb Q_p)$ that represents the element $\tilde w \in
\tilde W$ obtained as the image of $\nu$ under our injection $s:\Lambda_G
\to \tilde W$. [Since $\sigma(b)=b$ and $\tilde w$ has finite order modulo
$A_G(\mathbb Q_p)$, the element $b$ is basic. Moreover it is clear that
$\nu_b =\nu$.] 

So we may work with this particular $b$, the advantage of which is that 
$\Int(b)$ (conjugation by $b$) preserves both $A$ and the standard Iwahori
subgroup of $G$. Since $J_b$ is obtained from $G$ by twisting the Frobenius
automorphism $\sigma$ by $\Int(b)$, the twist $A_{\tilde w}$ of $A$ by
$\tilde w$ is an unramified maximal torus of $J_b$, whose split component
is a maximal split torus in $J_b$ (by Bruhat-Tits theory). Therefore
$\rk_{\mathbb Q_p}(J_b)$ is equal to the dimension of the split component
of $A_{\tilde w}$. Since
$\tilde w$ acts on $X_*(A)$ via
$w_\nu=p(\tilde w)$, we conclude that
\[
\rk_{\mathbb Q_p}(J_b)=\dim \mathfrak a^{w_\nu},
\]
as desired.

\bibliographystyle{amsalpha}


\providecommand{\bysame}{\leavevmode\hbox to3em{\hrulefill}\thinspace}

\end{document}